\documentclass [12pt]  {article} 
\usepackage{amssymb}
\def \CC{{\mathbb{C}}}

\def \FFF{{\cal{F}}}

\begin{document}

\begin{center}
{\Large {\bf Entire functions sharing simple $a$-points with their 
first derivative II}}\\
\bigskip
{\sc Andreas Schweizer\footnote{This paper was written while the 
author was supported by ASARC (South Korea).}}\\
\bigskip
{\small {\rm Department of Mathematics,\\
Korea Advanced Institute of Science and Technology (KAIST),\\ 
Daejeon 305-701\\
South Korea\\
e-mail: schweizer@kaist.ac.kr}}
\end{center}
\begin{abstract}
\noindent
We discuss some results around the following question: 
Let $f$ be a nonconstant complex entire function and 
$a$, $b$ two distinct complex numbers. If $f$ and its 
derivative $f'$ share their simple $a$-points and also 
share the value $b$, does this imply $f\equiv f'$?
\end{abstract}

\subsection*{1. Introduction}

Two meromorphic functions $f$ and $g$ are said to share the value 
$a\in\CC$ IM (ignoring multiplicity), or just to share the value 
$a$, if $f$ takes the value $a$ at exactly the same points as $g$,
or in other words, if the zeroes of $f(z)-a$ coincide with the zeroes 
of $g(z)-a$.
\par
If moreover for every such common zero $z_0$ the multiplicity of 
$z_0$ as a zero of $f(z)-a$ is the same as the multiplicity of $z_0$ 
as a zero of $g(z)-a$, then $f$ and $g$ are said to share the value 
$a$ CM (counting multiplicity).
\\ \\
{\bf Definition.}
Two meromorphic functions $f$ and $g$ 
{\bf share their simple $a$-points} 
if the points where $f$ takes the value $a$ with multiplicity one 
are exactly those where $g$ takes the value $a$ with multiplicity 
one, or in other words, if the simple zeroes of $f(z)-a$ coincide 
with the simple zeroes of $g(z)-a$.
\\ \\
Note that this does not imply any sharing of the multiple $a$-points.
So sharing simple $a$-points is weaker than sharing the value $a$ CM, 
but it cannot be directly compared to sharing the value $a$ IM.
\par
The most famous theorem in the context of entire functions sharing 
values with their derivative is
\\ \\
{\bf Theorem 1.1.} [MuSt, Satz 1] \it
Let $f$ be a nonconstant entire function and let $a$ and $b$ be complex 
numbers with $a\neq b$. If $f$ and $f'$ share the values $a$ and $b$ IM, 
then $f\equiv f'$.
\rm
\\ \\
See [YaYi, Theorem 8.3] for the same proof in English. 
See also [L\"uXuYi] for some later generalizations.
\par
In [Sch], adopting the approach from [L\"uXuYi], we proved some partial
results in the spirit of potential analogues of Theorem 1.1 for sharing 
of simple $a$-points.
\rm
\\ \\
{\bf Theorem 1.2.} [Sch] \it
Let $f$ be a nonconstant entire function and $a$, $b$ two distinct
complex numbers. Assume that $f$ and $f'$ share their simple $a$-points 
and also share their simple $b$-points.
\begin{itemize}
\item[(a)] If $ab=0$, then $f\equiv f'$.
\item[(b)] If $a+b=0$, then there are other solutions besides 
$f\equiv f'$, for example $f(z)=a\sin z$.
\item[(c)] If $ab(a+b)\neq 0$, it is not known whether $f\equiv f'$
necessarily holds or whether there are other solutions.
\end{itemize}
\rm
\bigskip
\noindent 
To conclude $f\equiv f'$, in general we need one more value.
\\ \\
{\bf Theorem 1.3.} [Sch, Theorem 2.3] \it
If a nonconstant entire function $f$ and its derivative $f'$ share
their simple $a_j$-points for three different complex numbers $a_1$,
$a_2$, $a_3$, then $f\equiv f'$.
\rm
\\ \\
Actually, in [Sch] we prove a version [Sch, Theorem 2.4] that 
is more general and slightly more technical than Theorem 1.3.
Similarly, the following is a slightly more general version of 
Theorem 1.2 (a).
\\ \\
{\bf Theorem 1.4.} [Sch, Theorem 2.2] \it
Let $f(z)$ be a nonconstant entire function and $0\neq a\in\CC$. 
If $f$ and $f'$ share their simple zeroes and if every simple 
$a$-point of $f$ is a (not necessarily simple) $a$-point of $f'$,
then $f\equiv f'$. 
\rm
\\ \\
Comparing part (b) of Theorem 1.2 with Theorem 1.1, sharing simple
$a$-points seems to be somewhat weaker than sharing $a$ IM.
\par
So we will investigate a mixture of the two:
For one value the simple points are shared, the other value is 
shared IM.
\\

\subsection*{2. The main results}

\noindent
If one of the two values is zero, we get a nice definitive result.
\\ \\
{\bf Theorem 2.1.} \it
Let $f$ be a nonconstant entire function and $0\neq a\in\CC$.
\begin{itemize}
\item[(a)] If $f$ and $f'$ share their simple zeroes and share the 
value $a$ IM, then  $f\equiv f'$.
\item[(b)] If $f$ and $f'$ share the value $0$ IM and share their 
simple $a$-points, then $f\equiv f'$.
\end{itemize}
\rm
\bigskip
\noindent
Completing the picture for the general case means giving an answer
(positive or negative) to the following
\\ \\
{\bf Question.}
Let $a$, $b$ be two distinct nonzero complex numbers. 
If a nonconstant entire function $f$ and its derivative $f'$ share 
their simple $a$-points and also share the value $b$ IM, does this 
imply $f\equiv f'$?
\\ \\
Currently the answer doesn't seem to be known. We cannot even offer
a guess. So we know even less than in Theorem 1.2, where at least
some counterexamples for the specific situation $a+b=0$ are known.
\par
The first part of the strategy in [L\"uXuYi] works nicely and gives 
the following intermediate result.
\\ \\
{\bf Proposition 2.2.} \it
Let $f$ be a nonconstant entire function and let $a$ and $b$ be two
distinct nonzero complex numbers. If $f$ and $f'$ share their simple 
$a$-points and also share the value $b$ IM, then $f$ is a transcendental
function of order at most one.
\rm
\\ \\
Under some stronger conditions we can prove what we want.
\\ \\
{\bf Theorem 2.3.} \it
Let $f$ be a nonconstant entire function and let $a$ and $b$ be two
distinct complex numbers. If $f$ and $f'$ share their simple 
$a$-points and also share the value $b$ CM, then $f\equiv f'$.
\rm
\\ \\
{\bf Example 2.4.}
A quick word on functions $f$ that are meromorphic in the whole complex 
plane. In that setting it is known ([Gu] and [MuSt2]) that if $f$ and 
$f'$ share two finite values CM, then $f\equiv f'$.
\par
On the other hand, the meromorphic function  
$$f(z)=z+2+\frac{1}{4z}$$ 
and its derivative $f'(z)=1-\frac{1}{4z^2}$ share their simple $1$-points 
and share the value $2$ CM. This shows that it is not possible to relax 
the sharing condition for one of the values from CM to sharing simple 
$a$-points. So Theorem 2.3 does not hold for meromorphic functions.
\\ \\
For the next result we need standard terminology from Nevanlinna theory.
In particular, $N_{1)}(r,\frac{1}{f-a})$, $N_{(2}(r,\frac{1}{f-a})$ and 
$\overline{N}_{(2}(r,\frac{1}{f-a})$ denote respectively 
the integrated counting function of the simple $a$-points of $f$, 
the integrated counting function of the multiple $a$-points of $f$,
and the integrated truncated counting function of the multiple 
$a$-points of $f$ where each such point is counted only once.
So $N_{(2}(r,\frac{1}{f-a})=N(r,\frac{1}{f-a})-N_{1)}(r,\frac{1}{f-a})$
and $\overline{N}_{(2}(r,\frac{1}{f-a})=
\overline{N}(r,\frac{1}{f-a})-N_{1)}(r,\frac{1}{f-a})$.
\\ \\
{\bf Proposition 2.5.} \it
Let $f$ be a nonconstant entire function and let $a$ and $b$ be two
distinct nonzero complex numbers. Assume that $f$ and $f'$ share their 
simple $a$-points and also share the value $b$ IM. If $f\not\equiv f'$,
then 
\begin{itemize}
\item[(a)] $T(r,f)\leq 2N_{(2}(r,\frac{1}{f-a})
+\overline{N}_{(2}(r,\frac{1}{f-a})+S(r,f)$;
\item[(b)] $N_{1)}(r,\frac{1}{f-a})+\frac{1}{2}N_{(2}(\frac{1}{f'-a})\leq 
N_{(2}(r,\frac{1}{f-a})+S(r,f)$.
\end{itemize}
\rm
\bigskip
\noindent
So if there is a $\lambda<\frac{2}{5}$ such that 
$N_{(2}(r,\frac{1}{f-a})\leq \lambda T(r,f)+S(r,f)$,
we can conclude $f\equiv f'$.
\\ \\
Coming back to Theorem 2.1:
That $f$ and $f'$ share their simple zeroes is obviously equivalent
to $f'$ having no simple zeroes and all zeroes of $f$ having at least 
multiplicity $3$. So the following result generalizes Theorem 1.4.
\\ \\
{\bf Theorem 2.6.} \it
Let $f$ be a nonconstant entire function and $0\neq a\in\CC$.
If every zero of $f$ has at least multiplicity $3$ and if every simple
$a$-point of $f$ is a (not necessarily simple) $a$-point of $f'$, then
$$f\equiv f'$$ 
or 
$$f(z)=\left(\frac{4a}{27\beta^2}e^{\frac{2}{9}z}+\beta e^{-\frac{1}{9}z}\right)^3$$
with a nonzero constant $\beta$.
\par
Conversely, every function of this form has the stronger property that
every simple $a$-point of $f$ is a simple $a$-point of $f'$.
\rm
\\ \\
From this we obtain another generalization of Theorem 1.2 (a).
\\ \\
{\bf Corollary 2.7.} \it
If a nonconstant entire function $f$ and its derivative $f'$ share their
simple $a$-points for some nonzero $a\in\CC$ and if every zero of $f$ has
at least multiplicity $3$, then $f\equiv f'$.
\rm
\\ \\
{\bf Example 2.8.}
The function  
$$f(z)=\frac{a}{2}(\sin(2z)+1)$$
shows that in Corollary 2.7 the condition on the zeroes cannot be
relaxed to multiplicity at least $2$.
\\ \\
The following results show that some sharing conditions can force 
$f\not\equiv f'$.
\\ \\
{\bf Theorem 2.9.} \it
Let $f$ be a nonconstant entire function and let $a$ and $b$ be two
distinct nonzero complex numbers. If every simple $a$-point of $f$ is 
a (not necessarily simple) $a$-point of $f'$, and if $b$ is a Picard
value of $f$, then
$$f(z)=Ce^{\frac{a}{a-b}z}+b$$
with a nonzero constant $C$.
\rm
\\ \\
{\bf Corollary 2.10.} \it
If a nonconstant entire function $f$ and its derivative $f'$ share their
simple $a$-points for some nonzero $a\in\CC$ and if $f$ has a finite Picard
value $b$ different from $a$, then
$$f(z)=Ce^{\frac{a}{a-b}z}+b$$
with a nonzero constant $C$.
\rm
\\

\subsection*{3. A normality criterion}

In order to prove Proposition 2.2 (and Corollary 2.3) we have 
to generalize a normality criterion from [L\"uXuYi].
\\ \\
{\bf Theorem 3.1.} \it
Let $\FFF$ be a family of holomorphic functions in a domain 
$D$. Let $a$ and $b$ be two finite complex numbers such that
$b\neq a, 0$, and let $K$ be a closed bounded subset of $\CC$
with $b\notin K$. 
\par
If, for each $f\in\FFF$ and $z_0 \in D$, 
$f(z_0)=a\Rightarrow f'(z_0)\in K$, and 
$f(z_0)=b\Rightarrow f'(z_0)=b$, 
then $\FFF$ is normal in $D$.
\rm
\\ \\
{\bf Proof.} \rm 
With the stronger condition $f=a\Rightarrow f'=a$ this is 
Theorem 1.3 from [L\"uXiYi], and our proof follows their proof 
closely. So we will be a bit sketchy and mainly emphasize
the points where one has to be careful with the more general 
conditions.
\par
We can assume that $D$ is the unit disk. If $\FFF$ is not normal, 
the family $\FFF_1 =\{f-a~:~f\in\FFF\}$ is not normal. Then, by 
the famous Zalcman Lemma [Za], there exist
\begin{itemize}
\item[(a)] a number $0<r<1$,
\item[(b)] points $z_n$ with $|z_n|<r$,
\item[(c)] functions $f_n\in\FFF$,
\item[(d)] positive numbers $a_n\to 0$,
\end{itemize}
such that
$$f_n(z_n +a_n\xi)-a=g_n(\xi)\to g(\xi)$$
locally uniformly, where $g(\xi)$ is a nonconstant entire function
on $\CC$.
\par
As in [L\"uXuYi] we see that the zeroes of $g(\xi)$ and the zeroes 
of $g(\xi)-(b-a)$ are all multiple. The only thing that is needed
for this is that $f'(z_0)$ is bounded for all $z_0 \in D$ with 
$f(z_0)=a$ resp. with $f(z_0)=b$. (Compare also the proof of 
Lemma 3.2 in [Sch].)
\par
Now comes the key point, namely showing that $g(\xi)$ never takes the 
value $b-a$. The proof of this is verbatim the same as in [L\"uXuYi], 
as it only requires the condition $f=b\Rightarrow f'=b$.
\par
So we have deduced that $g$ is a nonconstant entire function with 
one finite Picard value ($b-a)$ and one totally ramified value ($0$).
As such functions do not exist, $\FFF$ must have been normal.
\hfill$\Box$
\\ \\
{\bf Corollary 3.2.} \it 
Let $a$, $b$ be two complex numbers with $b\neq 0,a$.
Let $f$ be a nonconstant entire function with 
$f=a\Rightarrow f'\in\{a,0\}$ and $f=b\Rightarrow f'=b$.
Then $f$ is a transcendental function of order at most one.
\rm
\\ \\
{\bf Proof.} \rm 
Taking $K=\{a,0\}$ we see from Theorem 3.1 that the family
$\FFF=\{f_\omega(z)\ :\ \omega\in\CC\}$ with 
$f_\omega(z)=f(z+\omega)$ is normal on $\CC$. 
By [Mi, p.198 and p.211] this implies that 
the order of $f$ is at most $1$.
\par
Obviously a nonconstant polynomial $f$ cannot satisfy the condition
$f=b\Rightarrow f'=b$.
\hfill $\Box$
\\

\subsection*{4. The classical approach}

\noindent
{\bf Lemma 4.1.} \it 
Let $f$ be a nonconstant entire function with $f\not\equiv f'$ and let 
$a$ and $b$ be two distinct nonzero complex numbers. If $f$ and $f'$ 
share their simple $a$-points and also share the value $b$ IM, then 
$$T(r,f')\leq N_{(2}(r,\frac{1}{f-a})
+\overline{N}_{(2}(\frac{1}{f'-a})+S(r,f).$$
\rm
\bigskip
\noindent
{\bf Proof.} \rm 
Following the proof for the case $ab\neq 0$ in [MuSt] or 
[YaYi, Theorem 8.2] we get
\begin{eqnarray*}
N(r,\frac{1}{f-a})+N(r,\frac{1}{f-b}) & \leq & 
N(r,\frac{1}{f-f'})+N_{(2}(r,\frac{1}{f-a})\\
 & \leq & T(r,f)+N_{(2}(r,\frac{1}{f-a})+S(r,f).
\end{eqnarray*}
Adding
$$m(r,\frac{1}{f-a})+m(r,\frac{1}{f-b})\leq m(r,\frac{1}{f'})+S(r,f)$$
we obtain
$$T(r,f)\leq N_{(2}(r,\frac{1}{f-a})+m(r,\frac{1}{f'})+S(r,f).
\eqno{(*)}$$
On the other hand
\begin{eqnarray*}
N(r,\frac{1}{f'-a})+N(r,\frac{1}{f'-b}) & \leq &
N(r,\frac{1}{f-f'})+N(r,\frac{1}{f''})+\overline{N}_{(2}(r,\frac{1}{f'-a})\\
 & \leq & T(r,f)+N(r,\frac{1}{f''})+\overline{N}_{(2}(r,\frac{1}{f'-a})+S(r,f).
\end{eqnarray*}
Adding
$$m(r,\frac{1}{f'-a})+m(r,\frac{1}{f'-b})+m(r,\frac{1}{f'})\leq 
m(r,\frac{1}{f''})+S(r,f)$$
gives
$$2T(r,f')+m(r,\frac{1}{f'})\leq 
T(r,f)+T(r,f'')+\overline{N}_{(2}(r,\frac{1}{f'-a})+S(r,f),$$
which implies
$$T(r,f')+m(r,\frac{1}{f'})\leq 
T(r,f)+\overline{N}_{(2}(r,\frac{1}{f'-a})+S(r,f).$$
Adding the equation $(*)$ to this yields the desired inequality.
\hfill $\Box$
\\

\subsection*{5. Proofs of the main results}

\noindent
{\bf Proof of Theorem 2.1.} \rm \\ 
(a) This is just a weaker version of Theorem 1.4.\\
(b) Because of $f'=0\Rightarrow f=0$ we actually have
$f=a\Rightarrow f' =a$. By [L\"uXuYi, Corollary 1.1] this 
(together with $f=0\Rightarrow f'=0$) implies $f\equiv f'$ 
or $f=a(\frac{A}{2}e^{\frac{z}{4}}+1)^2$ for some nonzero constant 
$A$. But the second possibility is easily excluded, as then $f'$ 
has simple $a$-points (for $\frac{A}{2}e^{\frac{z}{4}}=1$) which are 
not $a$-points of $f$.
\hfill $\Box$
\\ \\
{\bf Proof of Proposition 2.2.} \rm \\ 
This follows immediately from Corollary 3.2.
\hfill $\Box$
\\ \\
{\bf Proof of Theorem 2.3.} \rm \\ 
If $a=0$ or $b=0$ we are done by Theorem 2.1. If not, we know from
Proposition 2.2 that the order of $f$ is at most one. But the famous
Br\"uck conjecture [Br] is proved for functions of finite order
[GuYa, Theorem 1], which means
$$\frac{f'-b}{f-b}\equiv c,$$
a nonzero constant. Solving this differential equation, we get
$f=\kappa e^{cz} +b-\frac{b}{c}$, where $\kappa$ is a nonzero constant. 
So $f'=\kappa ce^{cz}$ has simple $a$-points. Because of the sharing the 
differential equation then gives $c=1$.
\hfill $\Box$
\\ \\
{\bf Proof of Proposition 2.5.} \rm \\ 
(a) By the First Main Theorem we have
$$m(r,\frac{1}{f'})+N(r,\frac{1}{f'})\leq T(r,f')+S(r,f)$$
and 
$$N_{1)}(r,\frac{1}{f'-a})+N_{(2}(r,\frac{1}{f'-a})\leq T(r,f')+S(r,f).$$
Using
$$m(r,\frac{1}{f-a})\leq m(r,\frac{1}{f'})+S(r,f)$$
and 
$$N_{(2}(r,\frac{1}{f-a})-\overline{N}_{(2}(r,\frac{1}{f-a})
\leq N(r,\frac{1}{f'})$$
in the first equation, and
$$N_{1)}(r,\frac{1}{f'-a})=N_{1)}(r,\frac{1}{f-a})$$
and
$$2\overline{N}_{(2}(r,\frac{1}{f'-a})\leq N_{(2}(r,\frac{1}{f'-a})$$
in the second, and then adding them, we obtain
\begin{eqnarray*}
 &  &
m(r,\frac{1}{f-a})+N_{(2}(r,\frac{1}{f-a})
-\overline{N}_{(2}(r,\frac{1}{f-a})+N_{1)}(r,\frac{1}{f-a})
+2\overline{N}_{(2}(r,\frac{1}{f'-a})\\
 & \leq & 2T(r,f')+S(r,f).
\end{eqnarray*}
In combination with Lemma 4.1 this gives
$$T(r,f)-\overline{N}_{(2}(r,\frac{1}{f-a})
\leq 2N_{(2}(r,\frac{1}{f-a})+S(r,f).$$
This proves part (a).\\
(b) Using the sharing of simple $a$-points, the First Main Theorem and
Lemma 4.1, we get
$$N_{1)}(r,\frac{1}{f-a})+\frac{1}{2}N_{(2}(\frac{1}{f'-a})\leq 
N_{1)}(r,\frac{1}{f'-a})+N_{(2}(\frac{1}{f'-a})
-\overline{N}_{(2}(\frac{1}{f'-a})$$
$$\leq T(f')-\overline{N}_{(2}(\frac{1}{f'-a})+S(r,f)\leq 
N_{(2}(\frac{1}{f-a})+S(r,f).$$
\hfill $\Box$
\\ \\
{\bf Proof of Theorem 2.6.} \rm \\ 
We have 
$$\frac{(f')^2(f-f')}{f^2(f-a)}\equiv\gamma$$
for some constant $\gamma$.
This is shown in [Sch] in the proof of [Sch, Theorem 2.2], which
is our Theorem 2.6 with the additional condition that all zeroes 
of $f'$ are multiple. Up to that point, this additional condition 
has not been used in that proof.
\par
If $\gamma =0$ we have $f\equiv f'$. Otherwise we see from this
differential equation that all zeroes of $f$ have multiplicity
exactly $3$. So $f$ is a third power of another entire function
and, more conveniently, we can write 
$$f(z)=g^3(\frac{z}{3}).$$
Plugging this into the differential equation and cancelling some 
terms, we get
$$\frac{(g')^2(g-g')}{g^3-a}=\gamma.$$
By (the proof of) [Sch, Theorem 2.4] this implies
$$g(z)=\frac{4a}{27\beta^2}e^{\frac{2}{3}z}+\beta e^{-\frac{1}{3}z},$$
which proves the first claim. Writing
$$f=\frac{2^6 a^3}{3^9}u^2 
+\frac{2^4 a^2}{3^5}u+\frac{2^2 a}{3^2}+\frac{1}{u}$$
with $u=\frac{e^{\frac{1}{3}z}}{\beta^3}$, we have
$$f'=\frac{2^7 a^3}{3^{10}}u^2 +\frac{2^4 a^2}{3^6}u-\frac{1}{3u}$$
and
$$f''=\frac{2^8 a^3}{3^{11}}u^2 +\frac{2^4 a^2}{3^7}u+\frac{1}{9u}.$$
From this we see that
$f=a\Leftrightarrow u=\frac{27}{8a}$ or $u=\frac{-27}{a}$.
Moreover, $u=\frac{27}{8a}$ gives $f'=0$ whereas $u=\frac{-27}{a}$ 
gives $f'=a$ and $f''=\frac{23}{27}a\neq 0$.
This proves the extra claim.
\hfill $\Box$
\\ \\
{\bf Proof of Corollary 2.7.} \rm \\ 
Keeping the notation from the proof of Theorem 2.6 we see that
$f'=a$ if and only if $u=\frac{-27}{a}$ or 
$u=\frac{27}{16a}(5\pm 3\sqrt{3})$. But for the last two values
of $u$ we get $f''\neq 0$ and $f\neq a$.
So not every simple $a$-point of $f'$ is an $a$-point of $f$.
\hfill $\Box$
\\ \\
{\bf Proof of Theorem 2.9.} \rm \\ 
We trivially have $f=b\Rightarrow f'=b$. So by Corollary 3.2 the 
order of $f$ is at most $1$. But the only nonconstant entire functions 
of order at most one with Picard value $b$ are $f(z)=Ce^cz +b$ with 
nonzero constants $C$ and $c$. Since $f$ has simple $a$-points, the 
sharing condition forces $c=\frac{a}{a-b}$.
\hfill $\Box$
\\ \\
{\bf Proof of Corollary 2.10.} \rm \\ 
Theorem 2.9 for $b\neq 0$, and Corollary 2.7 for $b=0$. 
\hfill $\Box$
\\

\subsection*{\hspace*{10.5em} References}
\begin{itemize}

\item[{[Br]}] R.~Br\"uck: \rm On entire functions which share 
one value CM with their first derivative, \it Results in Math.
\bf 30 \rm (1996), 21-24

%
\item[{[Gu]}] G.~Gundersen: \rm Meromorphic functions that share
two finite values with their derivative, 
\it Pacific J. Math. \bf 105 \rm (1983), 299-309

\item[{[GuYa]}] G.~Gundersen and L.~Z.~Yang: \rm Entire functions 
that share one value with one or two of their derivatives, 
\it J. Math. Anal. Appl. \bf 223 \rm (1998), 88-95

\item[{[L\"uXuYi]}] F.~L\"u, J.~Xu and H.~Yi: \rm Uniqueness 
theorems and normal families of entire functions and their derivatives, 
\it Ann. Polon. Math. \bf 95.1 \rm (2009), 67-75

\item[{[Mi]}] D.~Minda: \rm Yosida functions, in: \it Lectures on
Complex Analysis, Xian 1987, (C.T. Chuang, ed.) \rm World Scientific,
Singapore, 1988, pp.197-213

\item[{[MuSt]}] E.~Mues and N.~Steinmetz: \rm Meromorphe Funktionen, 
die mit ihrer Ableitung Werte teilen, \it Manuscripta Math. \bf 29 \rm 
(1979), 195-206

\item[{[MuSt2]}] E.~Mues and N.~Steinmetz: \rm Meromorphe Funktionen, 
die mit ihrer Ableitung zwei Werte teilen, \it Resultate Math. \bf 6 \rm 
(1983), 48-55

\item[{[Sch]}] A.~Schweizer: \rm Entire functions sharing simple $a$-points 
with their first derivative, \it Houston J. Math. \bf 39 no. 4 \rm (2013), 
1137-1148 

\item[{[YaYi]}] C.-C.~Yang and H.-X.~Yi: \it Uniqueness theory
of meromorphic functions, \rm Kluwer Academic Publishers Group,
Dordrecht, 2003

\item[{[Za]}] L.~Zalcman: \rm A Heuristic Principle in Complex 
Function Theory, \it Amer. Math. Monthly \bf 82 \rm (1975), 813-817

\end{itemize}

\end{document}